%% file: mobius_examples.tex
\documentclass[preprint]{article}
\usepackage{amsthm}

\usepackage{epsfig}
\usepackage{amssymb}
\usepackage{graphicx}
\usepackage{fancyhdr}
\usepackage{amsmath}
\usepackage{mathrsfs}
\usepackage{algorithm}
\usepackage{algorithmic}

\numberwithin{equation}{section}

\newtheorem{deff}{Definition}[section]
\newtheorem{twr}{Theorem}[section]

\newtheorem{remark}{Remark}[section]

\newcommand{\comment}[1]{}

\newcommand{\I}{\mathbb I}
\newcommand{\e}[0]{\mathbf{e}}
\newcommand{\PP}[0]{\mathbf{P}}

\newcommand{\E}[0]{\mathbb{E}}

\newcommand{\C}[0]{\mathbf{C}}
\newcommand{\X}[0]{\mathbf{X}}

\newcommand{\ZZ}[0]{\mathbf{Z}}
\newcommand{\NB}[0]{\mathfrak{N}}

\newcommand{\monD}[0]{\downarrow}
\newcommand{\monU}[0]{\uparrow}
\newcommand{\trev}[0]{\overleftarrow}

\newcommand{\dual}[0]{*}

\begin{document}
\title{Strong Stationary Duality for M\"obius monotone Markov chains: examples}

 \author{Pawe{\l} Lorek\thanks{Work of both authors supported by NCN Research Grant DEC-2011/01/B/ST1/01305.} \footnote{\texttt{Email:} Pawel.Lorek@math.uni.wroc.pl}
\\{\it  University of Wroc{\l}aw} 
\and Ryszard Szekli\footnotemark[1]  \ \footnote{\texttt{Email:} Ryszard.Szekli@math.uni.wroc.pl} 
\\{\it  University of Wroc{\l}aw}  
}
%\date{\today}
\date{January 22, 2014}
\maketitle
\begin{abstract}
We construct strong stationary dual chains for Ising model on a circle, non-symmetric random walk on square lattice and a random walk on hypercube. The strong stationary dual chains are all sharp and have the same state space as original chains. We use M\"obius monotonicity of these chains with respect to natural orderings of the corresponding state spaces.
This method provides an easy way to find eigenvalues in the Ising model and  for a random walk on hypercube.
\medskip \par
\noindent
{\bf Keywords:} Markov chains; stochastic monotonicity; eigenvalues; M\"obius monotonicity; strong stationary duality; strong stationary times; separation distance; mixing time; Ising model; hypercube
\smallskip\par \noindent
{\bf AMS MSC 2010: } 60J10; 06A06; 60G40
\end{abstract}

%\tableofcontents
\section{Introduction}
\input{mob_introduction}

\section{M\"obius monotonicity and duality}\label{moebius}
\input{mob_monotonicity}

\section{M\"obius monotone Markov chains: examples}\label{moebius-ex}
\input{mob_ex_ising}

\input{mob_ex_weighted_graph}

\input{mob_ex_cube}

\section{Proofs}\label{proofs}
\input{mob_proof_ising}

\input{mob_proof_weighted_graph}

\input{mob_proof_cube}

% tymczasowo:
%\input{mob_inne}

\input{mob_biblio}
%%%%%%%%%%%%%%%%%%%%%%%%%%%%%%%%%%%%%%%%%%%%%%%%%%%%%%%%%%%%%%%%%%%%%%%%%%%%%%%%%%
 
\end{document}

%% file: mob_introduction.tex
%Introduction

Consider an ergodic Markov chain $\X=(X_n)_{n\ge 0}$ on a discrete (finite or countable) state space $\E$ with transition matrix $\PP$ and initial distribution $\nu$.
One way of studying the speed of convergence of $\X$ to its stationary distribution $\pi$ is to find (and bound its tail) so-called \textsl{Strong Stationary Time} (SST),
i.e. such a stopping time $T$ ($T$ implicitly depends on $\nu$) that it is independent from $X_T$, and $X_T$ has distribution $\pi$.
SST's  were introduced by Aldous and Diaconis \cite{AldDia86,AldDia87}, who also gave examples of SST and their applications. Many examples can also be found in Diaconis \cite{Diaconis_book}.
%(tu podac wiecej przykladow...) \par 
First examples of SST's were created by  \textsl{ad hoc} methods. A general approach was invented by  Diaconis and Fill \cite{DiaFil90} who introduced {\it dual processes}. They 
showed that for $\X$ there always exists absorbing, so-called \textsl{Strong Stationary Dual} (SSD) absorbing chain $\X^*$, such that its time  to absorption $T^*$ is equal, in distribution, to  a SST $T$ for $\X$.
Their proof is an existence type argument which  does not show how to construct a dual chain in general. 
They showed one tractable  case \cite[Theorem 4.6]{DiaFil90}, where the state space is linearly ordered. Under the 
condition of stochastic monotonicity (related to the linear order)   of the corresponding time-reversed chain (and some assumptions on the initial distribution) they gave a recipe of how to construct
a dual chain on the same state space. A special, and important, case is a stochastically monotone birth-and-death chain for which the dual chain is an absorbing birth-and-death chain. 
\par 
Strong stationary dual chains have a variety of applications. Diaconis and Fill \cite{DiaFil_ex} gave an extension of this theory to countable state spaces.
Fill \cite{Filla} gave a stochastic proof of a well-known theorem (usually attributed to Keilson),
which states that the first passage time from 0 to $M$ of a stochastically monotone birth-and-death process on $\{0,\ldots, M\}$ is equal, in distribution, to a sum of geometric random variables related to the spectral values of $\X$.
Similar results for continuous time birth-and-death processes were obtained by Diaconis and Miclo \cite{DiaMiclo}.
Diaconis and Saloff-Coste \cite{DiaSal} studied cut-off phenomena for birth-and-death chains using SSD theory. 
\par 
All the mentioned examples above (although very interesting) somehow rely on Theorem 4.6 of \cite{DiaFil90} which involves linearly ordered states space.
That is why most of the known examples are related to birth-and-death chains. The main underlying assumption is (classical)  stochastic monotonicity of the time-reversed chain.
Although this monotonicity is defined also for partially ordered state spaces, it is not sufficient  for an analogous  construction of a SSD  chain  as in Diaconis and Fill \cite{DiaFil90}.
Lorek and Szekli \cite{LorSze12}  gave a recipe of how to construct dual chains on partially ordered state spaces with a special feature that the duals have the same state space as original chains.
The assumption of the classical stochastic monotonicity was replaced by the assumption of M\"obius monotonicity. 
This extension (to partially ordered state spaces) opens a new way of finding SSD chains defined for not linearly ordered state spaces. 
The purpose of this paper is to get a new SSD insight to some classical examples of finite state Markov chains. 
In section \ref{moebius} we recall needed definitions and facts about M\"obius monotone chains. In section \ref{moebius-ex} we
present strong stationary duals for an Ising model on a circle, non-symmetric random walk on a square lattice, and for a random walk on a hypercube. 
We find eigenvalues in the 
Ising model and  for a random walk on hypercube immediately from the form of the duals - since they are pure birth chains.
We note in passing, that for Ising model we will always obtain pure birth dual chain for arbitrary graph.
In section \ref{proofs} we give proofs of the main results. We believe that the presented method should be applicable for many other examples and can be used to find bounds on the speed of convergence to stationarity, and to find cut-off phenomena.

%% file: mob_monotonicity.tex
In this section we  recall needed  results on SSD  and  M\"obius monotone chains.
For a more complete material on duality see Diaconis and Fill \cite{DiaFil90}, and for results on M\"obius monotone chains, see Lorek and Szekli \cite{LorSze12}.
%%%%%%%%%%%%%%%%%%%%%%%%%%%%%%%%%%%%%%%%%%%%%%%%%%%%%%%%%%%%%%%%%%%
\subsection{Strong Stationary Duality}\label{subsec:duality}
For an ergodic Markov chain $\X=(X_n)_{n\ge 0}$  with the transition matrix $\PP$ and initial distribution $\nu$,
we are  interested in bounding a distance between
$\nu\PP^k$ (a distribution of a chain at step $k$) and its stationary distribution $\pi$. 
Often used  distance  is the  total variation distance $d_{TV}(\nu\PP^k,\pi)=\max_{A\subset \E}|\nu\PP^k(A)-\pi(A)|$.
Another useful distance is the {\it separation distance} $s$  defined as follows: $s(\nu\PP^k,\pi)=\max_{\e\in\E} (1-\nu\PP^k(\e)/\pi(\e))$. For random times $T$ which are SST,  Aldous and Diaconis \cite{AldDia87} show that 
$d_{TV}(\nu\PP^k,\pi)\leq s(\nu\PP^k,\pi)\leq P(T>n)$. 
\par 
Let $\X^*$ be a Markov chain with transition matrix $\PP^*$, initial distribution $\nu^*$ and a state space $\E^*$,
with an absorbing state $\e_a^*$. 
Let $\Lambda\equiv\Lambda(\e^\dual,\e), \e^\dual\in \E^\dual, \e\in \E$ be a stochastic kernel (called   {\em a link}), such that  $\Lambda(\e_a^\dual,\cdot)=\pi$, for $\e^\dual_a\in \E^\dual$.
$\X^*$ is a Strong Stationary Dual  (SSD) chain for $\X$ if
\begin{equation}\label{eq:duality}
 \nu=\nu^\dual\Lambda \quad \mbox{ and } \quad \Lambda\PP=\PP^\dual\Lambda.
\end{equation}
Diaconis and Fill \cite{DiaFil90} proved that  the absorption time $T^*$ of $\X^*$ is a SST for $\X$.
Thus, the problem of  finding SST for $\X^*$ translates into the  problem of studying the absorption time of $\X^*$.
\par\noindent
\begin{deff}
Strong Stationary Dual chain $\X^*$ is called \textsl{sharp} if $s(\nu\PP^n,\pi)=P(T^*>n).$ 
\end{deff}

\begin{remark}\label{rem:eig}\rm
The relation (\ref{eq:duality}) implies that for finite $\E$ and $\E^\dual$, $\PP$ and $\PP^*$ have the same set of eigenvalues.
\end{remark}
It turns out, that in some examples we can easily identify the eigenvalues of $\PP^*$, and thus, by the above remark,
we will also obtain the eigenvalues of $\PP$ which are usually not easy to obtain directly. 
\medskip\par 
%\subsubsection{TO DO: czas do pochloniecia, gdy dane wsie wart. wlasne}
% \texttt{Moze udaloby sie takie cos: w Miclo ``On absorption time and Dirichlet eigenvalues''
% jest wzor na czas do pochloniecia jesli dane sa wszystkie wartosci wlasne - moze
% w dwoch przykladach, gdzie mamy te wart. wlasne udaloby sie wszystko policzyc?
% (najpierw trzeba tam to zrozumiec, to bedzie mieszanka geometrycznych)  }

%{\footnotesize
% We will also make use of the following lemma, usually attributed to Keilson.
% \begin{lem}
 % Let $\X$ be a Markov chain on finite state space $E=\{\e_1,\ldots,\e_N\}$ and transitions $\PP$.
 % Assume that
 % \begin{itemize}
  % \item 
  % there exist a particular absorption state $\e_a\in\E$, i.e., $\PP(\e_a,\e_a)=1$ such
 % that $\{\e_a\}$ is the unique irreducible class of $\PP^*$ (i.e., for any $\e\in \E$ there exists
 % a path going from $\e$ to $\e_a$, $\e=\e_0,\e_1,\ldots,\e_q=\e_a$, satisfying $\PP(\e_i,\e_{i+1})>0$ for 
 % $i\in\{0,\ldots,q-1\}$. Moreover, assume that submarkoviav restriction $\bar\PP:=(\PP(\e,\e'))_{\e,\e'\in \bar\E})$
 % of $\PP$ to $\bar\E:=\E\setminus\{\e_a\}$ is irreducible (i.e., any pair of points from $\bar\E$ 
 % can be joined by a path).
%  
 % \item $\PP$ is reversible w.r.t some positive probability meausure on $\bar\E$ 
 % \item All eigenvalues of $\PP$ are nonnegative, i.e., they fulfill
 % $1>\lambda_0>\lambda_1\geq \lambda_2\geq\ldots\geq \lambda_{N-1}>0$
 % \end{itemize}
% \end{lem}
%}

%%%%%%%%%%%%%%%%%%%%%%%%%%%%%%%%%%%%%%%%%%%%%%%%%%%%%%%%%%%%%%%%%%%

\subsection{Duality for M\"obius monotone chains}\label{subsec:mobius}
%Note that in Section \ref{subsec:duality} no recipe for dual was given.
%We can make a use of a duality \textsl{once} $\Lambda, \nu^*, \E^*$ and $\PP^*$ are found. 
In this section we recall how to construct a SSD  chain for finite partially ordered state spaces. 
We shall consider a finite state space $\E=\{\e_1,\ldots,\e_M\}$ with a partial ordering $\preceq$. From the very beginning
we shall choose an enumeration of $\E$ such that $\e_i\preceq\e_j$ implies $i<j$ (which is always  possible).  We call such an enumeration consistent with $\preceq$.
With this enumeration the partial  ordering can by represented by an upper-triangular, 0-1 valued matrix $\C$.
The inversion $\C^{-1}$ represents  (in the incidence algebra) the so  called \textsl{M\"obius function}, usually denoted by $\mu$, see Rota \cite{rota}.
The  M\"obius function allows for the following  calculus:  it is possible to  recover $f$ from the relation $\bar{F}(\e)=\sum_{\e:\e\succeq \e_i} f(\e)$, namely
$f(\e_i)=\sum_{\e:\e\succeq \e_i} \mu(\e_i,\e)\bar{F}(\e)$.
%\mymarginpar{\tiny chyba lepiej zrezygnowac z dwoch mon. (strzalki do gory i na dol), a skupic sie na jednej}

\begin{deff}
Let $\PP$ be a transition matrix with enumaration of states consistent with $\C$. We say that $\PP$ (or alternatively, $\X$) is 
${}^\monD$-M\"obius monotone (\ ${}^\monU$-M\"obius monotone) if $\C^{-1}\PP\C \geq 0$ (\ $(\C^T)^{-1}\PP\C^T \geq 0$) (each entry is nonnegative).
\end{deff}
\noindent  We say that  $\mathbf{f}:\E\to\mathbb{R}$ is ${}^\monD$-M\"obius monotone (${}^\monU$-M\"obius monotone)
if $\mathbf{f} (\C^{T})^{-1}\geq 0$ ($\mathbf{f} \C^{-1}\geq 0$).
In terms of the transition probabilities, we have 
\par \medskip
\begin{tabular}{ccc}
${}^\monD$-M\"obius monotonicity: & $\forall(\e_i,\e_j\in\E)\quad \sum_{\e:\e\succeq \e_i} \mu(\e_i,\e)\ \PP(\e,\{\e_j\}^\downarrow)\geq 0$,\\ 
\ & \ \\
${}^\monU$-M\"obius monotonicity: & $\forall(\e_i, \e_j\in\E)\quad \sum_{\e: \e\preceq \e_j} \PP(\e,\{\e_i\}^\monU)\mu(\e,\e_j) \geq 0$,\\
 \end{tabular}
 %\mymarginpar{\tiny czy w tych sumach \par dawac $P$ czy $\PP$?}
%\par \medskip
\par \medskip 
\noindent 
where $\{\e_j\}^\monD=\{\e:\e\preceq \e_j\}$,  $\{\e_j\}^\monU=\{\e:\e\succeq \e_j\}$, and  $\PP(\e,A)=\sum_{\e'\in A} \PP(\e,\e')$.

\par \medskip 
\noindent 
We recall the SSD  result of Lorek and Szekli \cite{LorSze12} ($\trev\X$ denotes the time-reversed process). 
%%+++++++++++++++++++++++++++++++++++++++++++++++++++++++++
%%++++++++++++++++++++++++++++++++++++++++++++++++++++++
\begin{twr}[Lorek and Szekli \cite{LorSze12}]\label{lorsze:mobius_main_twr}
Let $\X$ be an ergodic Markov chain on a finite state space $\E=\{\e_1,\ldots,\e_M\}$, which is partially ordered with $\preceq$, and  has a unique maximal state $\e_M$. For   the stationary distribution $\pi$ and an initial distribution $\nu$ we
assume that 
 \begin{itemize}
 \item[(i)] $g(\e)={\nu(\e)\over \pi(\e)}$  is $^\monD$-M\"obius monotone,
 \item[(ii)] $\trev\X$ is $^\monD$-M\"obius monotone.\label{eq:rev_Mob}
\end{itemize}
Then there exists a Strong Stationary Dual chain $\X^\dual$ on $\E^\dual=\E$ with link being a truncated stationary distribution
$\Lambda(\e_j,\e_i)=\I(\e_i\preceq \e_j) {\pi(\e_i)\over H(\e_j)},$
where $H(\e_j)=\sum_{\e:\e\preceq \e_j}\pi(\e)$. The initial distribution and transitions of $\X^*$ are given, respectively, by
\begin{eqnarray}
 \nu^\dual(\e_i)&=&H(\e_i)\sum_{\e:\e\succeq \e_i}\mu(\e_i,\e)g(\e),\nonumber\\
\PP^\dual(\e_i,\e_j)& =&
 \frac{H(\e_j)}{ H(\e_i)}\sum_{\e:\e\succeq \e_j} \mu(\e_j,\e)\trev \PP(\e,\{\e_i\}^\monD).\label{eq:dual}
\end{eqnarray}
\end{twr}
%+++++++++++++++++++++++++++++++++++++++++++++++++++++++++++
\begin{remark}\label{rem:dual_sharp}\rm Following Remark 2.39 of Diaconis and Fill \cite{DiaFil90} and the terminology used there,
the  Strong Stationary Dual  $\X^\dual$ in Theorem \ref{lorsze:mobius_main_twr} is sharp, and the corresponding strong stationary time is the time to stationarity, i.e., $s(\nu\PP^n,\pi)=P(T>n).$ The reason for this is that $\Lambda(\e^\dual,\e _M)=0$ for all, $\e^\dual \neq \e _M \in \E^\dual$.
\end{remark}
\begin{remark}\rm
 Theorem \ref{lorsze:mobius_main_twr} is stated for $^\monD$-M\"obius monotonicity, but it can be similarly stated for $^\monU$-M\"obius monotonicity
 (see Corollary 3.1 in \cite{LorSze12}).  The other formulation is potentially useful, because a chain can be, e.g., $^\monD$-M\"obius monotone but not $^\monU$-M\"obius montone.
\end{remark}
\begin{remark}\rm
The assumption on the initial distribution is not very restrictive, for example  if $\e_1$ is a unique minimal state and $\nu=\delta_{\e_1}(\cdot)$, then the assumption is fulfilled, and also $\nu^*=\delta_{\e_1}(\cdot)$. 
For simplicity of presentation, in all subsequent examples the initial distribution will be the single  atom at the minimal element (this assumption may be relaxed).
\end{remark}
%The assertion of Theorem \ref{lorsze:mobius_main_twr} can be given in a matrix form: 
 %\begin{equation}\label{dual:matrix_form}
%\PP^*=\diag(\pi\C)^{-1} (\C^T\diag(\pi)) \PP (\C^T\diag(\pi))^{-1} \diag(\pi\C), 
 %\end{equation}
 %where $\diag(\mathbf{v})$ denotes a diagonal matrix with diagonal entries $v_1,\ldots,v_M$, %$\mathbf{v}=(v_1,\ldots,v_M)$.
 In order to find and use the above constructed SSD chains one has   to find \textsl{an approperiate} ordering (w.r.t which the chain is M\"obius monotone). It is  worth mentioning, that for linearly ordered state space $^\monU$-M\"obius monotonicity is equivalent to the usual stochastic monotonicity,
  in general partially ordered spaces this is not the case.  It
 turns out that  for partially ordered spaces some \textsl{natural} orderings work. A non-symmetric  random walk on the unit cube is an example presented in \cite{LorSze12}. In the next section we shall give new examples.
 %Moreover, above matrix form easily allows for  checking the monotonicity on computer (for $$\E$ of reasonable size): simply check if for given ordering (i.e., given $\C$) the matrix
% $\PP^*$ in (\ref{dual:matrix_form}) has nonnegative entries. 
 %
 %(tu dokladniej: ze to pierwszy krok: jesli nie to nie jest monotoniczny, jesli tak,
 %to duza szansa, ze jest)

%% file: mob_ex_ising.tex
\subsection{Ising model on a circle}
Let $G=(V,E)$ be a finite graph. Elements of state space $\E=\{-1,1\}^V$ are called \textsl{configurations},
and for $\e\in\E$ the value $\e(v)$ is called the \textsl{spin} at vertex $v$.   
For a given configuration $\e$ its \textsl{energy} is defined as
$$\mathcal{H}(\e)\quad=-\sum_{\{x,y\}\in E} \e(x)\cdot\e(y),$$
where the sum is over all edges of the graph. For $\beta\geq 0$, the \textbf{Ising model} on the graph $G$ with
parameter $\beta$ is the probability measure on $\E$ given by
\begin{equation}\label{eq:pi_ising}
\pi(\e)={ e^{-\beta \mathcal{H}(\e)} \over Z_\beta},
\end{equation}
where $Z_\beta=\sum_{\e\in\E} e^{-\beta \mathcal{H}(\e)}$ is a  normalizing constant. The parameter $\beta$
has a physical interpretation as the inverse of the  temperature of the configuration. Note, that for $\beta=0$ (equivalent to infinite temperature), every spin configuration is equally likely,
i.e., it is the same as setting spin at each vertex to -1 or +1 with probability  1/2 independently. 
In general,  $\beta$ represents the influence of energy $\mathcal{H}$ on $\pi$.
\par

\par 
This model has  focused a  lot of attention in the context of speed of convergence to equilibrium of particle systems.
Propp and Wilson \cite{ProWil} introduced \textsl{Coupling From The Past} algorithm and used it 
to show how  to draw exact sample from (\ref{eq:pi_ising}) in the case of square lattice.
%\mymarginpar{\tiny   moze te opisy o \par   ``lots of attention''\par    dac do introduction?}
Recently Ding and Peres \cite{DinPer} showed that for   Ising models on each graph it takes at least $(1/4+o(1))n\log n$ steps for the Glauber dynamics  to mix,
where $n$ is the corresponding number of vertices. In Ding and Peres \cite{DinPer2} a simple proof for the bound $n\log n/2 $ was presented.

We shall consider the   Ising model  on a circle.
We  set $V=\{0,\ldots,N-1\}$ and $E=\{\{i,(i+1)\mod N\}: i=0,\ldots,N-1\}$.
%(everytnig would work in a similar way for Ising model on a line with edges $E=\{(i,i+1): i=0,\ldots,N-1\}$).
%From now on, we always consider vertices number modulo $N$ and vertex $0-1$ means $N$.
%\mymarginpar{\tiny  nie brzmi dobrze:\par ``means $N$''...}
The  distribution (\ref{eq:pi_ising}) in this case can be written as
\begin{equation}\label{eq:pi_ising_line}
\pi(\e)={1\over Z_\beta}\exp\left(\beta\sum_{i=0}^{N-1} \e(i)\e(i+1)\right). 
\end{equation}
(we always mean vertex number modulo $N$).
We build a \textsl{Gibbs sampler} for this model with stationary distribution (\ref{eq:pi_ising_line}). 
This chain has the state space $\E=\{-1,1\}^V$ and its dynamics can be described as follows
\begin{itemize}
 \item Given a configuration $\e$ at step $n$, i.e., $X_n=\e$, choose a vertex $i\in\{0,\ldots,N-1\}$ with probability $1/N$.
\item  Take $U_{n+1}$,  a random variable with the uniform distribution $U(0,1)$, independent of $U_i,  i\leq n$.
Update the spin at vertex $k$ in the following way
$$X_{n+1}(i)=\left\{
\begin{array}{lll}
 +1 & \mathrm{ if } \ U_{n+1}< \displaystyle {e^{2\beta(k_+(i,\e)-k_-(i,\e))}\over e^{2\beta(k_+(i,\e)-k_-(i,\e))}+1},\\[6pt]
-1 & \mathrm{ otherwise, }
\end{array}\right.
$$
where $k_+(i,\e)$ is the number of neighbours of vertex $i$, in configuration $\e$, with spin values +1, and $k_-(i,\e)$ is the number of neighbours
of vertex $i$, in configuration $\e$, with spin values -1. 
\end{itemize}
Note that, for the  circle, we have $k_+(i,\e), k_-(i,\e) \in\{0,1,2\}$, and $k_+(i,\e)+k_-(i,\e)=2$.
The chain $\X$ constructed in this way is   reversible. Moreover, $\X$  can be viewed as a random walk on $N-$dimensional cube, where the probability 
of changing coordinate $i$ depends  on the values of the neighbouring coordinates.
 
It turns out that if we consider the coordinate-wise ordering,  i.e.,  
$\e\preceq \e'$ if $\e(i)\leq \e'(i)$ for every vertex $i\in V$, then $\X$ is $^\monD$-M\"obius monotone.  Let $M:=2^{|V|}=2^N$.
Denote by $\e_1$ the state with all spins equal to $-1$ (minimal state), and by $\e_{M}$ the state with all spins equal to $+1$
(maximal state). We identify $\E=\{-1,1\}^V$ with the enumerated set $\{\e_1,\ldots,\e_M\}$, where the enumeration is consistent with $\preceq$.
% It will be convienient to define $$\e_k^i:=(\e(0),\ldots,\e({k-1}),i,\e({k+1}),\ldots,\e({N-1})),\  i\in\{-1,+1\}.$$
Applying Theorem \ref{lorsze:mobius_main_twr} we obtain
\begin{twr}\label{twr:Ising}
Consider the Gibbs sampler $\X$ for the Ising model on the circle with vertices $V=\{0,\ldots,N-1\}$, and edges 
$E=\{\{i,i+1 \mod N\}: i=0,\ldots,N-1\}$. Assume that $\X$ starts with the configuration $\e _1$. 
Then, there exists  sharp SSD chain $\X^*=(X_n^*)_{n\ge 0}$ on the state space $\E^*=\E$, with the state $\e_{M}$ being the absorbing one, 
 starting with probability 1 at $\e_1$, and having  transition probabilities for $\e,\e'\in\{\e_1,\ldots,\e_M\}$
\begin{equation}\label{Dual_ising}
\PP^*(\e,\e')=
 \left\{
\begin{array}{lll}
 0 & \mathrm{if} & \e\succ \e'\\[8pt]
 \displaystyle  {1\over N} S(\e) & \mathrm{if} & \e=\e'\\[8pt]
 \displaystyle {H(\e')\over H(\e)} {1\over N}\left(1-\displaystyle {e^{2\beta(k_+(j,\e)-k_-(j,\e))}\over e^{2\beta(k_+(j,\e)-k_-(j,\e))}+1}\right)   & \mathrm{if} & \e'=\e+s_j,\ \e(j)=-1\\
\end{array}
\right.
\end{equation}
where $s_j=(0,\ldots,0,2,0,\ldots,0)$ (2 on the $j$-th coordinate), $S(\e)=\sum_{i=0}^{N-1} \pmb 1 \{\e(i)=1\}$, and $H(\e)=\sum_{\e'\preceq \e} \pi(\e)$.
\end{twr}
%+++++++++++++++++++++++++++++++++++++++++++++++++++++++
Note  that our SSD chain $\X^*$ jumps, with probability 1,  only to greater or equal states in the ordering $\preceq$, thus its eigenvalues are the entries on the diagonal of the matrix $\PP^*$ written using the enumeration of the states consistent with this ordering.  Therefore the eigenvalues are $\{0,1/N,2/N,\ldots,1\}$, and from Remark
\ref{rem:eig} these values are also the eigenvalues of the original  transition matrix $\PP$ for the Gibbs sampler $\X$. The multiplicity of eigenvalue $k/N$ is $N \choose k$.
\par 
Consider one dimensional projection $Z^*_n:=S(X^*_n)$. Note that $\ZZ^*=(Z^*_n)_{n\ge 0}$ is also a Markov chain, since for all $j=1\ldots,2^N$, we have
$$\PP^*(\e,[\e_j])=\PP^*(\e',[\e_j])\qquad \textrm{for all\ } \e\sim_S\e',$$
where $\e\sim_S\e'$ iff $S(\e)=S(\e')$, $\e,\e'\in \E$,
and $[\e_j]:=\{\e: S(\e)=S(\e_j)\},\ j=0,\ldots,N-1$ denote the equivalence classes of the relation $\sim_S$. 
\par 
%\mymarginpar{\tiny mozna bez tych klas rownowaznosci, od razu pisac jak dziala}
With $X^*_0=\e_1$, and thus with $Z^*_0=0$, the time to absorption $T^*$ for $\X^*$ is the same
as the time to absorption at $N$ for $\ZZ^*$. Note that $\ZZ^*$ is a pure-birth chain on the state space $\{0,\ldots, N\}$, with birth rates $\lambda_k=1-{k\over N}$, $k=0,\ldots,N-1$, and absorption at $N$.
The time to absorption $T^*$ is thus equal, in distribution, to $\sum_{i=0}^{N-1} Y_i,$ where $Y_i, \ i=0,\ldots, N-1$ are independent random variables, and $Y_i$ has
geometric distribution (on $\{1,\ldots \}$) with the success parameter $p_i=1-{i\over N}$. A  coupon-collector argument shows that, for $n=N\log N+cN$, we have
 (equality  because of Remark \ref{rem:dual_sharp})
$$s(\delta_{\e_1}\PP^n,\pi)=P(T^*>n)\leq e^{-c},\  c>0 .$$

%% file: mob_ex_weighted_graph.tex
\subsection{Random walk on weighted directed graph}
Consider a random walk on a  directed weighted graph $G=(V,E)$ with vertices $V=\{v_1,v_2,\ldots,v_n\}$, edges $E=\{(i,j): $ edge from $v_i$ to $v_j\}$ and
with a weighting function $w: E\to [0,\infty)$. Denote by $w_{i,j}$ the nonnegative weight of the directed edge
from  node $v_i$ to $v_j$. If there is no edge between these nodes, i.e., $(i,j)\notin E$, then $w_{i,j}=0$. 
We allow $w_{i,i}$ be  nonzero.
\par 
Let $\mathcal{N}(i)=\{j: (i,j)\in E\}$ be a set of neighbours of node $v_i$. 
Random walk may be viewed as a process of sequential vertex visiting. 
We assume that weights are normalized, i.e., for all $i\in\{1,\ldots,n\}$ we have $w_{i,i}+ \sum_{r\in\mathcal{N}(i)} w_{i,r}=1$.
The  probability of a single step from node $i$ to $j$ is then given by  $P(i,j)= w_{i,j}$.
\par 
In this section we consider the following example: Let $V=\{0,1,\ldots,N\}^2$ with edges 
\begin{equation}\label{rw:graph_E}
 % ((x_1,y_1),(x_2,y_2))\in E \iff ( |x_1-x_2|=1 \textrm{ and } y_1=y_2) \textrm{ or } (x_1=x_2 \textrm{ and } |y_1-y_2|=1)
 ((x_1,y_1),(x_2,y_2))\in E \iff  |x_1-x_2|+ |y_1-y_2|=1
\end{equation}
for  $ x_1,x_2,y_1,y_2\in\{0,\ldots,N\}$. 
Thus, for each node there are at most four edges in four directions: \textsl{up, down, left, right} plus a possible self-loop.
The weighting function depends only on the direction in the following way: for $((x_1,y_1),(x_2,y_2))\in E$ and nonnegative parameters 
$\lambda_1, \lambda_2, \mu_1, \mu_2$
such that $\lambda_1 + \lambda_2+ \mu_1+ \mu_2 \leq 1$
\begin{equation}\label{rw:graph_weights}
w_{((x_1,y_1),(x_2,y_2))}=
\left\{
\begin{array}{lll}
 \lambda_1 & \mathrm{if} & x_2=x_1+1, y_2=y_1, \\[3pt]
 \mu_1 & \mathrm{if} & x_2=x_1-1, y_2=y_1, \\[3pt]
 \lambda_2 & \mathrm{if} & x_2=x_1, y_2=y_1+1, \\[3pt]
 \mu_2 & \mathrm{if} & x_2=x_1, y_2=y_1-1, \\ [3pt]
 \displaystyle 1-\sum_{(x,y)\in \mathcal{N}((x_1,y_1))} w_{((x_1,y_1),(x,y))} & \mathrm{if} & x_2=x_1, y_2=y_1. \\ 
\end{array}
\right.
\end{equation}
We associate weights directly with one step probabilities: $$\PP((x_1,y_1),(x_2,y_2))=w_{((x_1,y_1),(x_2,y_2))}.$$
%Without loss of generality we can assume that $\lambda_1+\lambda_2+\mu_1+\mu_2=1$, then the weights can be directly associated with probabilities of moving along the edge.
Roughly speaking, we consider a random walk on square lattice $\{0,\ldots,N\}^2$, at each step we can move  (if feasible): \textsl{right} with probability $\lambda_1$,
\textsl{left} with probability $\mu_1$, \textsl{up} with probability $\lambda_2$ and \textsl{down} with probability $\mu_2$. With remnant probability we stay at a given vertex. For convenience, we let $\rho_1:=\lambda_1/\mu_1$, and $\rho_2:=\lambda_2/\mu_2$.
Denote the transition matrix of a corresponding Markov chain $\X$ by $\PP$. 
The chain is time-reversible (i.e. $\trev{\PP}=\PP$) and has (time-reversibility equations can be easily checked) the stationary  distribution  on $V$
$$\pi((x,y))=C^{-1} \rho_1^x \rho_2^y$$
for $(x,y)\in V=\{0,\ldots,N\}^2$,
where the normalizing constant  $C$ for $\rho_1\neq 1$ and $\rho_2\neq 1$ is given by $$C={ 1-\rho_1^{N+1}\over 1-\rho_1} \cdot { 1-\rho_2^{N+1}\over 1-\rho_2},
$$
and $C$ for other cases can be obtained by obvious modifications.
 
We shall use the  coordinate-wise partial ordering $(x_1,y_1)\preceq (x_2,y_2) \iff x_1\leq x_2 \textrm{ and } y_1\leq y_2.$
Then we have  unique minimal element $\e_1=(0,0)$  and the  maximal one $\e_M=(N,N)$, where $M=(N+1)^2$.
It turns out that $\X$ is  M\"obius monotone  for any set of parameters $\lambda_1,\mu_1,  \lambda_2,\mu_2 >0$, such that $\lambda_1+\lambda_2+\mu_1+\mu_2 \le 1$, and applying Theorem \ref{lorsze:mobius_main_twr} we have:
%%%+++++++++++++++++++++++++++++++++++++++++++++++++++++
%%+++++++++++++++++++++++++++++++++++++++++++++++++++++++++
\begin{twr}\label{twr:graph}
Let $\X$ be a random walk on directed weighted graph with $G=(V,E)$, with $V=\{0,\ldots,N\}^2$, and $E$ given in (\ref{rw:graph_E}), weights given in (\ref{rw:graph_weights}) and
with positive parameters $\lambda_1\neq \mu_1, \ \lambda_2\neq \mu_2$, such that $\lambda_1+\lambda_2+\mu_1+\mu_2\le 1$. Assume, that $\X$ starts at $\e_1=(0,0)$.
Then there exists sharp SSD   chain $\X^*$ which is an absorbing Markov chain (with $\e_M=(N,N)$ being the single absorbing state) on the state space $\E^\dual =\E=\{0,\ldots,N\}^2$, starting at $\e_1=(0,0)$, with 
the following transition probabilities (for $x, x', y, y' \in \{0,\ldots,N\}$)
\begin{equation}\label{eq:dual_graph}
\PP^*((x,y),(x',y'))=
\left\{\begin{array}{lll}
\frac{1-\rho_1^{x+2}}{1-\rho_1^{x+1}}\cdot  \mu_1 & \mathrm{ if\ } x'=x+1,\  y'=y\\[8pt]
\frac{1-\rho_2^{y+2}}{1-\rho_2^{y+1}}\cdot  \mu_2 & \mathrm{ if\ } y'=y+1,\  x'=x\\[8pt]
\frac{1-\rho_2^{y}}{1-\rho_2^{y+1}}\cdot  \lambda_2 & \mathrm{ if\ } x'=x, \ y'=y-1, \ y'\neq N \\[8pt]
\frac{1-\rho_1^x}{1-\rho_1^{x+1}}\cdot  \lambda_1 & \mathrm{ if\ } y'=y, \ x'=x-1, \ x'\neq N \\[8pt]
 1-(\lambda_1+\lambda_2+\mu_1+\mu_2) & \mathrm{ if\ }  x'=x, y'=y, \ (x,y)\in \{0,\ldots ,N-1\}^2 \\[8pt]
1-(\lambda_2+\mu_2)& \mathrm{if}  x'=x, y'=y,\ y=N, \ x\in\{0,\ldots ,N-1\}\\[8pt]
1-(\lambda_1+\mu_1)&\mathrm{if}  x'=x, y'=y,\ x=N, \ y\in\{0,\ldots ,N-1\}\\[8pt]
 1 	& \mathrm{ if\ } x'=x=y=y=N \\[8pt]
\end{array}\right.
\end{equation}
\end{twr}
%%%++++++++++++++++++++++++++++++++++++++++++++++++++++++++
Thus, the SSD chain $\X^*$ is again a chain on $\E$, with feasible moves in the same directions  as $\X$ \textsl{except} for movements on the upper borders of this square lattice.
Once the chain hits the border $(\cdot,N)$ ( or $(N,\cdot)$), then it can only move \textsl{left} or \textsl{right} (\textsl{up} or \textsl{down}) until it hits the absorbing state $(N,N)$. 
Note  that probability of changing $i$-th coordinate, $i=1,2$, is independent of the value of $(3-i)$-th coordinate. The chain $\X^*$, for a suitable selection of the parameters, can have a {\it drift} towards the absorbing state. Note that the case $\rho_1=1$, and/or $\rho_2=1$ can be obtained by obvious modifications in computing $H(x,y)$ (see the proof in section \ref{proof1}).

%\mymarginpar{\bf \tiny example  with drift???}

One can study the  time to absorption $T^*$ in the following way: it is the time of hitting a border $(\cdot,N)$ or $(N,\cdot)$ plus the time for the one dimensional birth-and-death chain with 
birth probability  $\lambda_1$ and death probability $\mu_1$ (or $\lambda_2$ and $\mu_2$ respectively) to reach the state $N$ (worst cases scenarios can be used).

%%%%%%%%%%%%%%%%%%%%%%%%%%%%%%%%%%%%%%%%%%%%%%%%%%%%%%%%%%%%%%%%%%%%%%

%% file: mob_ex_cube.tex
\subsection{Random change of single  coordinate on  a cube}
Let us consider a discrete time Markov chain $\X$ with state space $ \E=\{0,\ldots,k\}^n$, 
which evolves in the following way: it stays with probability $1/2$ or  (with probability $1/2$) for   one coordinate chosen uniformly, it changes uniformly its value  to any other different value. In terms of the transition probabilities, for $\e=(\e (1),\ldots,\e (n))\in \E, \ \e(i) \in\{0,\ldots,k\}$, we set
\begin{equation}\label{eq:dna_chain}
\PP(\e, \e')=
\left\{\begin{array}{llr}
{1\over 2} &\textrm{ if }   \e=\e'\\[5pt]
{1\over 2nk}  & \textrm{ if for some} \   i&  \e (i)\neq \e '(i) \ \textrm{and}\ \  \e (j)=\e '(j), \ j\neq i  \\[5pt]
0 & \textrm{otherwise}
\end{array}
\right. 
\end{equation}
\par 
Since $\PP$ is symmetric, the corresponding stationary distribution  is uniform, i.e.,  $$\pi(\e)={1\over (k+1)^n},\ \e\in\E.$$
The motivation for this example comes from DNA sequence alignment. Given $n$ sequences of length $k+1$ the task
is to find points of references in each one such that, starting \textsl{reading} sequence $i$ from it's reference point $r(i)$
we have the biggest agreement in all sequences. Since the state space is huge (of size $(k+1)^n$), often Monte Carlo methods are used.
One constructs a chain such that its stationary distribution assigns higher mass to states with high agreements.
The chain given in (\ref{eq:dna_chain}) is a simplified version of such a chain. 
\par 
The chain $\X$ can be seen as en extension of the standard lazy random  walk on the unit cube (obtained for $k=1$).
Using   the coordinate-wise ordering $\preceq$ on $\E$, it turns out that $\X$ (which is reversible) is M\"obius monotone.
%To present the dual chain let us introduce some notation first. 
For this  ordering, the state $\e_1=(0,\ldots,0)$ is the minimal state
and $\e_M=(k,\ldots,k)$ is the maximal state (with $M=(k+1)^n$), where we use an enumeration of $\E$ consistent with $\preceq$.
%Let $\e_i^{(k)}=(0,\ldots,0,k,0,\ldots,0)$ be a state with $k$ on position $i$, and zeros otherwise.
Applying Theorem \ref{lorsze:mobius_main_twr} we obtain.
%+++++++++++++++++++++++++++++++++++++++++++++++++++++++++
%+++++++++++++++++++++++++++++++++++++++++++++++++++++++++++
\begin{twr}\label{twr:cube}
 Consider the chain $\X$ described above, on state space $\E=\{0,\ldots,k\}^n$,   with transition probabilities  given in (\ref{eq:dna_chain}).
 Assume that $\X$ starts at $\e_1$.
 Then, there exists sharp SSD chain $\X^*$ on the  state space $\E^*=\E$,  with the state $\e_{M}$ being the absorbing one, 
 starting with probability 1 at $\e_1$, and having  transition probabilities, for all $A\subseteq\{1,\ldots,n\}$, $j\notin A$
 $$\begin{array}{lcl}
\PP^*(\e_A^{(k)},\e_{A\cup j}^{(k)}) & = & \displaystyle{(k+1)\over 2nk}, \\[10pt]
 \PP^*(\e_A^{(k)},\e_A^{(k)}) & = & \displaystyle{n(k-1)+|A|(k+1)\over 2nk}, 
  \end{array}
$$
 where $\e_A^{(k)}=(\e (1),\ldots,\e (n))$ with $\e(i)=k$ if $i\in A$ and $\e (i)=0$ if $i\notin A$, and all other transitions have probability 0.\par
\end{twr}
%++++++++++++++++++++++++++++++++++++++++++++++++++++++++++
Note  that SSD chain $\X^*$ jumps, with probability 1,  only to greater or equal states in the ordering $\preceq$, thus its eigenvalues are the entries on the diagonal of the matrix $\PP^*$ written using an enumeration of the states consistent with this ordering.
The
states which can be traversed by $\X^*$ are of  the form $\e_A^{(k)}$, which means that  $\X^*$ can be identified with  a random walk on the unit cube $\{0,k\}^n$.
Again, by Remark \ref{rem:eig}, the eigenvalues of $\PP$ are the same as diagonal entries of $\PP^*$,
i.e., 
$$ {n(k-1)+i(k+1)\over 2nk}, \quad i=0,1,\ldots,n.$$
%The dual chain starting at $\e_{min}$ can only be in states of form $\e_A^{(k)}$ for some $A\subseteq \{1\ldots,n\}$. It means, that
%actual state space of the dual is $\E^*=\{0,k\}^*$ and thus de facto it is a random walk on cube. \par 
Similarly as in the  Ising model example, we can consider the time to absorption of one dimensional projection $Z_t^*:=S(X^*_t),$ where
$S(\e)=\sum_{i=1}^n \pmb 1\{\e(i)=k\}$. If $Z^*_0=0$, then the time to absorption $T^*$ of $Z^*_t$ is the same
as for $X^*_t$, and is distributed as the sum of independent variables $\sum_{i=0}^{n-1} Y_i,$ where $Y_i$ has geometric distribution with the success parameter $p_i={(n-i)(k+1)\over 2nk}$.
For the expected absorption time we have
$$ET^*=\sum_{i=0}^{n-1} {1\over p_i}= \sum_{i=0}^{n-1}{1\over n-i}{2nk\over k+1} = {2nk\over k+1}\sum_{i=1}^n{1\over i} \leq {2k\over k+1} (n+1)\log n $$ 
For the variance of $T^*$ we have
$$VarT^*=\sum_{i=0}^{n-1} {1-p_i\over p_i^2}={2nk\over (k+1)^2}\sum_{i=0}^{n-1}{nk-n+ki+i\over (n-i)^2}
$$
$$={2nk\over (k+1)^2}\left[nk \sum_{i=0}^{n-1}{1\over (n-i)^2}+k\sum_{i=0}^{n-1}{i\over (n-i)^2}- \sum_{i=1}^n {1\over i}\right] \stackrel{(*)}{\leq} \left({2nk\over k+1}\right)^2 {\pi^2\over 6},$$
where in $^{(*)}$ we used the following ineqalities
$$\sum_{i=0}^{n-1}{1\over (n-i)^2} \leq  {\pi^2\over 6}, \quad \sum_{i=0}^{n-1}{i\over (n-i)^2} \leq n{\pi^2\over 6}.$$
By Remark \ref{rem:dual_sharp} and from Chebyshev's inequality, we have that after $m={2k\over k+1} (n+1)\log n + c {2k\over k+1} {\pi\over\sqrt{6}} n, c\geq 0 $ steps we have
%\marginpar{\tiny uwaga: moze lepiej dac: $m={2k\over k+1} (n+1)\log n + d n$?}\par
$$s(\nu\PP^m,\pi)=P(T>m) \leq  P(T-ET \leq c\sqrt{Var})\leq P(|T-ET| \leq c\sqrt{Var})\leq {1\over c^2}.$$

%% file: mob_proof_ising.tex
\subsection{Proof of Theorem \ref{twr:Ising}}
The M\"obius function for the  coordinate-wise ordering is given by
$$\mu(\e,\e')
=\left\{
\begin{array}{ll}
 (-1)^{S(\e')-S(\e)} & \mathrm{if \ } \e\preceq\e',\\[5pt]
 0 & \mathrm{otherwise},
\end{array}
\right.
$$
where $S(\e)=\sum_{i=0}^{N-1} \pmb{1}\{\e(i)=+1\}$. 
For convienience, let us define
$$f(i,\e):={ e^{2\beta(k_+(i,\e)-k_-(i,\e))}\over e^{2\beta(k_+(i,\e)-k_-(i,\e))}+1},$$
so that the probability of choosing vertex $i$ and updating it to +1 in configuration $\e$ is ${1\over N} f(i,\e)$.
 %$$ f(k_+({ e^{2\beta(a_+-a_-)}\over e^{2\beta(a_+-a_-)}+1}$$
 Let us define also $$\e_i^r:=(\e(0),\ldots,\e({i-1}),r,\e({i+1}),\ldots,\e({N-1})),\  r\in\{-1,+1\}$$
 and we will shortly write $\e_i^-$ for $\e_i^{-1}$ and $\e_i^+$ for $\e_i^{+1}$. Moreover, by $\NB(i)$ we denote the neighbours of vertex $i$, i.e., $\NB(i)=\{i-1,i+1\}$.
 We will directly apply Theorem \ref{lorsze:mobius_main_twr}. We can replace $\trev{\PP}$ by $\PP$, since this chain is reversible.
 \medskip\par \noindent
First, let us calculate $\PP^*(\e_i^+,\e_i^-)$.
$$
\begin{array}{rl}
  \displaystyle{H(\e_i^+)\over H(\e_i^-)}\PP^*(\e_i^+,\e_i^-)&  \displaystyle=\sum_{\e\succeq \e_i^-}\mu(\e_i^-,\e)\PP(\e,\{\e_i^+\}^\monD)\\[16pt]
  \displaystyle = \ \PP(\e_i^-,\{\e_i^+\}^\monD)  & -\displaystyle  \sum_{j: \e_i^-(j)=-1 \atop j\notin\NB(i),j\neq i} \PP(\e_i^-+s_j,\{\e_i^+\}^\monD)
- \sum_{j: \e_i^-(j)=-1 \atop j \in\NB(i)} \PP(\e_i^-+s_j,\{\e_i^+\}^\monD)  \\[16pt]
\displaystyle -\  \PP(\e_i^-+s_i,\{\e_i^+\}^\monD)  & +\displaystyle\sum_{j: \e_i^+(j)=-1 \atop j \notin\NB(i),j\neq i} \PP(\e_i^++s_j,\{\e_i^+\}^\monD)
+\sum_{j: \e_i^+(j)=-1 \atop j \in\NB(i)} \PP(\e_i^++s_j,\{\e_i^+\}^\monD) 
\end{array}
$$
Since $\e_i^-+s_i=\e_i^+$, we have
$${H(\e_i^+)\over H(\e_i^-)}\PP^*(\e_i^+,\e_i^-)= \PP(\e_i^-,\{\e_i^+\}^\monD) - \PP(\e_i^-+s_i,\{\e_i^+\}^\monD) $$
$$- \sum_{j: \e_i^-(j)=-1 \atop j\notin\NB(i),j\neq i} {1\over N} \left[1-f(j,\e_i^-+s_j)\right] 
- \sum_{j: \e_i^-(j)=-1 \atop j \in\NB(i)} {1\over N}\left[1-f(j,\e_i^-+s_j)\right]  $$
\begin{equation}\label{eq:Nb1}
+\sum_{j: \e_i^+(j)=-1 \atop j \notin\NB(i),j\neq i}{1\over N}\left[1-f(j,\e_i^++s_j)\right] 
+\sum_{j: \e_i^+(j)=-1 \atop j \in\NB(i)} {1\over N}\left[1-f(j,\e_i^++s_j)\right] . 
\end{equation}
Note that for $j\notin \NB(i)$ we have $k_+(j,\e_i^-+s_j)=k_+(j,\e_i^++s_j)$ and $k_-(j,\e_i^-+s_j)=k_-(j,\e_i^-+s_j)$, thus
$f(j,\e_i^-+s_j)=f(j,\e_i^++s_j)$ in this case. The corresponding terms with summation over $j\notin\NB(i)$ cancel out. For the first two terms we have
$$ \PP(\e_i^-,\{\e_i^+\}^\monD) = 
1- \sum_{j: \e_i^-(j)=-1\atop j\notin\NB(i),j\neq i} {1\over N} f(j,\e_i^-+s_j)
-\sum_{j: \e_i^-(j)=-1\atop j\in\NB(i)} {1\over N} f(j,\e_i^-+s_j)$$
and
$$ \PP(\e_i^+,\{\e_i^+\}^\monD)  = 
1- \sum_{j: \e_i^+(j)=-1\atop j\notin\NB(i)} {1\over N} f(j,\e_i^++s_j)
-\sum_{j: \e_i^-(j)=-1\atop j\in\NB(i)} {1\over N} f(j,\e_i^++s_j).$$
Again, for $j\notin \NB(i)$, we have $k_+(j,\e_i^-+s_j)=k_+(j,\e_i^++s_j)$, and $k_-(j,\e_i^-+s_j)=k_-(j,\e_i^-+s_j)$,
thus $f(j,\e_i^-+s_j)=f(j,\e_i^++s_j)$. The corresponding terms again cancel out in $ \PP(\e_i^-,\{\e_i^+\}^\monD) - \PP(\e_i^+,\{\e_i^+\}^\monD)$.
Plugging in the remaining sums to (\ref{eq:Nb1})   we obtain
$$
\begin{array}{lll}
 \displaystyle{H(\e_i^+)\over H(\e_i^-)}\PP^*(\e_i^+,\e_i^-) = & - &\displaystyle
  \sum_{j: \e_i^-(j)=-1\atop j\in\NB(i)} {1\over N} f(j,\e_i^-+s_j)
- \sum_{j: \e_i^-(j)=-1 \atop j \in\NB(i)} \left[1-f(j,\e_i^-+s_j)\right] \\[16pt]
 & + & \displaystyle \sum_{j: \e_i^-(j)=-1\atop j\in\NB(i)} {1\over N} f(j,\e_i^++s_j)+\sum_{j: \e_i^+(j)=-1 \atop j \in\NB(i)} \left[1-f(j,\e_i^++s_j)\right].  
\end{array}
$$
Now all the terms on the right hand side in the above equality  cancel out, therefore $\PP^*(\e_i^+,\e_i^-)=0$.
\medskip\par 
\noindent
Next, we get immediately
$$\PP^*(\e_i^-,\e_i^+)=
{H(\e_i^+)\over H(\e_i^-)}\sum_{\e\succeq \e_i^+}\mu(\e_i^+,\e)P(\e,\{\e_i^-\}^\monD)={H(\e_i^+)\over H(\e_i^-)}{1\over N} f(i,\e_i^-).
$$
 \par 
\noindent
We have yet to calculate  the probability of staying at $\e$.
$$
\begin{array}{rl}
\PP^*(\e,\e)=&\displaystyle
\sum_{\e'\succeq \e}\mu(\e,\e')\PP(\e',\{\e\}^\monD)=P(\e,\{\e\}^\monD)-\sum_{j: \e(j)=-1}\PP(\e+s_j,\e)\\[16pt]
=&\displaystyle 1 - \sum_{j: \e(j)=-1} \PP(\e,\e+s_j) -\sum_{j: \e(j)=-1}\PP(\e+s_j,\e)\\[16pt]
=&\displaystyle 1 - \sum_{j: \e(j)=-1} {1\over N} f(j,\e)  -\sum_{j: \e(j)=-1}{1\over N} \left[1-f(j,\e+s_j)\right]\\[16pt]
=&\displaystyle 1 - \sum_{j: \e(j)=-1} {1\over N} f(j,\e)  +\sum_{j: \e(j)=-1}{1\over N} f(j,\e+s_j) -\sum_{j: \e(j)=-1}{1\over N}\\[16pt]
=&\displaystyle 1- {N-S(\e)\over N}={S(\e)\over N}.
\end{array}
$$
Summing up, we obtain $\PP^*$ given in (\ref{Dual_ising}).

%% file: mob_proof_weighted_graph.tex
\subsection{Proof of Theorem  \ref{twr:graph}}\label{proof1}

We start with a detailed expression for the transition probabilities of $\X$ 
$$\PP((x,y),(x',y'))=
\left\{\begin{array}{ll}
 \lambda_1 & \mathrm{ if\ } x'=x+1\leq N, y'=y\\[5pt]
 \lambda_2 & \mathrm{ if\ } x'=x, y'=y+1\leq N\\[5pt]
\mu_1 & \mathrm{ if\ } x'=x-1\geq 0, y'=y\\[5pt]
\mu_2 & \mathrm{ if\ } y'=y-1\geq 0, x'=x\\[5pt]
1-(\lambda_1+\lambda_2+\mu_1+\mu_2) & \mathrm{if\ } x'=x>0, y'=y>0\\[5pt]
1-(\lambda_1+\lambda_2+\mu_1) & \mathrm{if\ } x'=x>0, y'=y=0\\[5pt]
1-(\lambda_1+\lambda_2+\mu_2) & \mathrm{if\ } x'=x=0, y'=y>0\\[5pt]
1-(\mu_1+\mu_2) & \mathrm{if\ } x'=x=y=y'=N\\
1-(\mu_1+\mu_2+\lambda_1) & \mathrm{if\ } x'=x>0, y'=y=N\\[5pt]
1-(\mu_1+\mu_2+\lambda_2) & \mathrm{if\ } x'=x=N, y'=y>0\\[5pt]
1-(\lambda_1+\lambda_2) & \mathrm{if\ } x'=x=y=y'=0\\
\end{array}\right.
$$
In a standard way we can check that $\X$ is reversible and the stationary distribution is given by
%{\bf Stationary distribution:} 
$$\pi((x,y))=C^{-1} \rho_1^x\rho_2^y$$
where  $C$ is the normalizing constant, and $\rho_i=\lambda_i/\mu_i, \ i=1,2$.
%Process is time-reversible, i.e. $\trev{\PP}=\PP $\par
For the coordinate-wise ordering
%{\bf Ordering:} 
$$(x,y)\preceq (x',y') \iff x\leq x' \textrm{ and } y\leq y' ,$$
with  the minimal state $\e_1=(0,0)$, and the maximal state $\e_M =(N,N)$, ($M=(N+1)^2$)
directly from Proposition 5 in \cite{rota}, we find the corresponding
M\"obius function: 
$$\begin{array}{llrll}
 \mu((x,y),(x,y)) & = & 1 \\[4pt]
 \mu((x,y),(x+1,y)) & = & -1 & x+1\leq N\\[4pt]
 \mu((x,y),(x,y+1)) & = & -1 & y+1\leq N\\[4pt]
 \mu((x,y),(x+1,y+1)) & = & 1 & x+1\leq N, y+1\leq N \\[4pt]
  & = & 0 & \mathrm{otherwise} .\\[4pt]
\end{array}
$$
For 
$$H(x,y)=C^{-1}\sum_{x'\leq x} \rho_1^{x'}  \sum_{y'\leq y} \rho_2^{y'}
= C^{-1} (1-\rho_1)^{-1} (1-\rho_2)^{-1} \left(1-\rho_1^{x+1}\right) (1-\rho_2^{y+1}),$$
we shall compute
\begin{equation}\label{dual-transitions}
\PP^*((x,y),(x_2,y_2))={H(x_2,y_2)\over H(x,y)} \sum_{(x',y')\succeq (x_2,y_2)} \mu((x_2,y_2),(x',y')) \trev{\PP}((x',y'),\{(x,y)\}^\downarrow).
\end{equation}
Set
 $$ S:= \sum_{(x',y')\succeq (x_2,y_2)} \mu((x_2,y_2),(x',y')) \trev{\PP}((x',y'),\{(x,y)\}^\downarrow).$$
Note that in order to prove that $\trev\X$ is $^\monD$-M\"obius monotone it is enough to show that $S\ge 0$. Since $\X$ is reversible, we take $\PP$ instead of $\trev{\PP}$ in the above formula. We shall consider all possible transitions, case by case.
\begin{itemize}
\item (inside lattice, up $x$ direction)
%Najpierw rozpatrzymy sytuacje poza brzegami:
 
$x_2=x+1,\   y_2=y$
 $$S=\sum_{(x',y')\succeq (x+1,y)} \mu((x+1,y),(x',y'))  {\PP}((x',y'),\{(x,y)\}^\downarrow)$$
% sumujemy po $(x',y')\succeq (x+1,y)$, 
where $\mu$ will be non-zero only in the following cases
 $$  \mu((x+1,y),(x+1,y))=1, \quad \mu((x+1,y),(x+1,y+1))=-1, $$
 $$\mu((x+1,y),(x+2,y))=-1, \quad \mu((x+1,y),(x+2,y+1))=1.$$  
Combining these cases with the values of   ${\PP}((x',y'),\{(x,y)\}^\downarrow)$ we get 
$$ S= \mu((x+1,y),(x+1,y))\PP((x+1,y),\{(x,y)\}^\downarrow)-1\cdot 0 -1\cdot 0 +1\cdot 0 = \mu_1,$$
 \item (inside lattice, up $y$ direction)

$x_2=x,\   y_2=y+1$

 in a similar way as ubove, we get
$$ S= \mu((x,y+1),(x,y+1))\PP((x,y+1),\{(x,y)\}^\downarrow)-1\cdot 0 -1\cdot 0 +1\cdot 0 = \mu_2,$$
\item (inside lattice, down $x$ direction)

$x_2=x-1\geq 0,\  y_2=y$
 
 %Tu wiecej opcji: najlepiej rysunek sobie zrobic: prosokat o rogach (0,0) i $(x,y)$ (to jest %$\{x,y\}^\downarrow$)
 %na nim punkt $(x-1,y)$ - z tego punktu ``startuje'' $\mu$ i sa cztery opcje (+0 lub +1 na %kazdej wspolrzednej)
using the formula for $S$ we have
$$ S= \mu((x-1,y),(x-1,y))\PP((x-1,y),\{(x,y)\}^\downarrow)
+\mu((x-1,y),(x,y))\PP((x,y),\{(x,y)\}^\downarrow)$$
$$
+\mu((x-1,y),(x,y+1))\PP((x,y+1),\{(x,y)\}^\downarrow)
\mu((x-1,y),(x,y+1))\PP((x,y+1),\{(x,y)\}^\downarrow)
$$
$$
=1\cdot (1-\lambda_2)-1\cdot(1-\lambda_2-\lambda_1)-1\cdot\mu_2+1\cdot\mu_2=\lambda_1,
$$
 \item (inside lattice, down $y$ direction)

$x_2=x,\  y_2=y-1\geq 0$
$$ S= \mu((x,y-1),(x,y-1))\PP((x,y-1),\{(x,y)\}^\downarrow)
+\mu((x,y-1),(x,y))\PP((x,y),\{(x,y)\}^\downarrow)$$
$$
+\mu((x,y-1),(x+1,y-1))\PP((x+1,y-1),\{(x,y)\}^\downarrow)
+ \mu((x,y-1),(x+1,y))\PP((x+1,y),\{(x,y)\}^\downarrow)
$$
$$
=1\cdot (1-\lambda_1)-1\cdot(1-\lambda_2-\lambda_1)-1\cdot\mu_1+1\cdot\mu_1=\lambda_2,
$$

\item (inside lattice, down on both axes)

$x_2=x-1\ge 0,\  y_2=y-1\ge 0$
$$ S= \mu((x-1,y-1),(x-1,y-1))\PP((x-1,y-1),\{(x,y)\}^\downarrow)
+\mu((x-1,y-1),(x-1,y))\PP((x-1,y),\{(x,y)\}^\downarrow)$$
$$+\mu((x-1,y-1),(x,y-1))\PP((x,y-1),\{(x,y)\}^\downarrow)
+
\mu((x-1,y-1),(x,y))\PP((x,y),\{(x,y)\}^\downarrow)
$$
$$
=1\cdot 1 -1\cdot (1-\lambda_2)-(1-\mu_1)+\mu_1+\mu_2=0.
$$
%np.
%\item $x_2=x-2, y_2=y-1$
%$$ S= \mu((x-2,y-1),(x-2,y-1))\trev{\PP}((x-2,y-1),\{x,y\}^\downarrow)$$
%$$
%+\mu((x-2,y-1),(x-1,y-1))\trev{\PP}((x-1,y-1),\{x,y\}^\downarrow)$$
%$$+\mu((x-2,y-1),(x-1,y-1))\trev{\PP}((x-1,y-1),\{x,y\}^\downarrow)
%+
%\mu((x-2,y-1),(x-2,y))\trev{\PP}((x-2,y),\{x,y\}^\downarrow)
%$$
%$$
%=1\cdot 1 -1\cdot 1 -1 \cdot (1-\lambda_2)+1\cdot (1-\lambda_2)=0
%$$
In a similar way it is possible to check that inside the lattice the only one remaining movement with positive probability is the feedback movement

\item (feedback inside lattice)

$x_2=x>0, \ y_2=y>0$

$$\PP^*((x,y),(x,y))=1-\lambda_1-\lambda_2-\mu_1-\mu_2 =\PP((x,y),(x,y)),$$

 \item (upper border, up $x$ direction)

$x_2=x+1\leq N,\ y_2=y=N$ 
$$S=\mu( (x+1,N), (x+1,N)) \PP((x+1,N),\{(x,N)\}^\downarrow)=\mu_1,$$

\item (upper border, down $y$ direction)

$x_2=x<N, y=N, y_2=N-1$ 
$$ S= \mu((x,N-1),(x,N-1))\PP((x,N),\{(x,y)\}^\downarrow)
+\mu((x,N-1),(x,N)\PP((x,y),\{(x,y)\}^\downarrow)$$
$$
+\mu((x,N-1),(x+1,N-1))\PP((x+1,N-1),\{(x,y)\}^\downarrow)
\mu((x,N-1),(x+1,N))\PP((x+1,N),\{(x,y)\}^\downarrow)
$$
$$
=1\cdot (1-\lambda_1) - 1\cdot(1-\lambda_1)-1\cdot \mu_1+1\cdot\mu_1=0,
$$

\item (upper border, down $x$ direction)

$x_2=x-1\geq 0, y_2=y=N$
$$ S=\mu((x-1,N),(x-1,N))\PP((x-1,N),\{(x,y)\}^\downarrow)+\mu((x-1,N),(x,N)) \PP((x,N),\{(x,y)\}^\downarrow)$$
$$=+1\cdot 1-1\cdot(1-\lambda_1)=\lambda_1,$$

\item (lower border,  up $x$ direction)

$x_2=x+1\leq N, y_2=y=0$ 
$$S=\mu((x+1,0),(x+1,0)  \PP((x+1,0),\{(x,0)\}^\downarrow)=\mu_1,$$

\item (lower border, down $x$ direction) 

$x_2=x-1\geq 0, y_2=0$ 
$$S=\mu((x,0),(x-1,0))=
 \mu((x-1,0),(x-1,0)) \PP((x-1,0),\{(x,0)\}^\downarrow)
+\mu((x-1,0),(x,0)) \PP((x,0),\{(x,0)\}^\downarrow)
$$
$$
+\mu((x-1,0),(x-1,1) \PP((x-1,1),\{(x,0)\}^\downarrow)
+\mu((x-1,0),(x,1)) \PP((x,1),\{(x,0)\}^\downarrow)
$$
$$
= 1\cdot(1-\lambda_2)-1\cdot(1-\lambda_1-\lambda_2)-1\cdot \mu_2+\mu_2=\lambda_1,$$

\item (lower border, up $y$ direction) 

$x_2=x\geq 0, y_2=1, y=0$ 
$$S=\mu((x,1),(x,1)) \PP((x,1),\{(x,0)\}^\downarrow)=\mu_2.$$
In a similar way we get
%++++++++++++++++++++++++++++\par 

 \item (right border, up $y$)

$x_2=x=N, y_2=y+1\leq N$,\  \ \  
$S\mu_2,$

\item (right border, down $y$)

$x_2=x=N, y_2=y-1\leq N$,\ \ \    
$S=\lambda_2,$

\item (right border, down $x$)

$x_2=N-1, x=N, y_2=y<N$,\ \ \   
$S=0,$

\item (left border, up $y$)

$x_2=x=0, y_2=y+1\leq N$,\  \ \   
$S=\mu_2,$

\item (left border, up $x$)

$x_2=x+1\leq N, y_2=y $,\ \ \   
$S=\mu_1,$

\item (left border, down $y$)

$x_2=x=0, y_2=y-1\geq N$,\ \ \    
$S= \lambda_2,$

\item (absorbing state)

$x_2=x=N, y_2=y=N$,\ \ \ 
$ S = 1.$

\item (feedback movements)

for all $(x,y)\in \{0,\ldots , N_1\}^2$,\ \ \ $S=1-(\lambda_1+\lambda_2+\mu_1+\mu_2,$

for $x=N$, and $y\in\{0,\ldots,N-1\}$,\ \ \ $S=1-(\lambda_2+\mu_2),$

for $y=N$, and $x\in\{0,\ldots , N-1\}$,\ \ \ $S=1-(\lambda_1+\mu_1).$
\end{itemize}

Now using (\ref{dual-transitions}), and using values of $H(x,y)$ we obtain $\PP^*$ given in (\ref{eq:dual_graph}). 
%Possible transitions for $\PP^*$ are illustrated in Figure \ref{fig:2stations_independent_dual}.

%\begin{figure}
%\centering
%\includegraphics[width=14cm]{pics/2stations_independent_dual.eps}
%\caption{$\PP^*$ transitions. }
%\label{fig:2stations_independent_dual}
%\end{figure}

%%%%%%%%%%%%%%%%%%%%%%%%%%%%%%%%%%%%%%%%%%%%%%%%%%%%%%%%%%%%%%%%%%%%%%%%%%%%%%%%%%

%% file: mob_proof_cube.tex
\subsection{Proof of Theorem \ref{twr:cube}}
Consider the coordinate-wise ordering
%\centerline{ $\e = (\e(1),\ldots, \e(n)) \preceq (\e'(1),\ldots,\e'(n))= \e'$ iff $\e(i)\leq \e'(i), i=1,\ldots,n$.}
$$\e = (\e(1),\ldots, \e(n)) \preceq (\e'(1),\ldots,\e'(n))= \e' \textrm{ iff }\e(i)\leq \e'(i), i=1,\ldots,n.$$

Again, for this ordering with minimal element $\e_1=(0,\ldots,0)$ and maximal element $\e_M=(k,\ldots,k)$ (with $M=(k+1)^n$) , directly from Proposition 5 in Rota \cite{rota}, we find the
corresponding M\"obius function
$$
\mu((\e(1),\ldots,\e(n)),(\e(1)+d_1,\ldots,\e(n)+d_n))
=\left\{\begin{array}{cr}
  \displaystyle(-1)^{\sum_{i=1}^n d_i}  &  d_i\in\{0,1\}, \ \e(i)+d_i\leq n, \\
 & i=1,\ldots,n  \\[8pt]
0 & \textrm{otherwise}. \\
\end{array}\right.
$$
For $H(\e)=\sum_{\e'\preceq \e} \pi(\e')= |\{\e': \e'\leq \e\}|\cdot 1/(k+1)^n$, we shall compute directly transitions of the dual chain (\ref{eq:dual}) from  Theorem \ref{lorsze:mobius_main_twr}.
Note, that in order to prove that $\trev{\X}$ is ${}^\monD$-M\"obius monotone, it is enough to show that all summands in (\ref{eq:dual}) are non-negative. We take $\PP$ instead of $\trev\PP$ 
since this chain is reversible.
%$$\PP^*(\e_i,\ej)={H(\e')\over H(\e)} \sum_{\e:\e\preceq \e'} \mu(\e',\e)\trev{\PP}(\e,\{\e\}^\monD)$$

\par 
\noindent For convenience, we shall consider states of the following form
 $$ \e_A^{(k)}=(\e_A^{(k)}(1),\ldots,\e_A^{(k)}(k)), \ \ A\subseteq \{1,\ldots,n\},$$ with $\e_A^{(k)}(i)=k$ if $i\in A$ and 0 otherwise.
Note, that there are $(k+1)^{|A|}$ states smaller or equal (w.r.t. $\preceq$) to $\e_A^{(k)}$, and we have
\begin{equation}\label{eq:HeAk}
 {H(\e_{A\cup\{j\}}^{(k)})\over H(\e_A^{(k)})} = {(k+1)^{|A\cup\{j\}|}\over (k+1)^{|A|}}=k+1\ \textrm{for } j\notin A.
\end{equation}
Let us calculate transitions of the dual chain from state $\e_A^{(k)}$. 
We shall use $s_i=(0,\ldots,0,1,0,\ldots,0)$ with $1$ at position $i$. For the probability of staying at this state we get
$$\PP^*(\e_A^{(k)},\e_A^{(k)})=1\cdot   \sum_{\e\succeq \e_{A}^{(k)}} 
\mu(\e_{A}^{(k)},\e)\PP(\e,\{\e_A^{(k)}\}^\monD)$$
$$=\mu(\e_{A}^{(k)},\e_{A}^{(k)})\PP(\e_A^{(k)},\{\e_A^{(k)}\}^\monD)
+\sum_{i\in \{A\}^c} \mu(\e_{A}^{(k)},\e_{A}^{(k)}+s_i) \PP(\e_{A}^{(k)}+s_i,\{\e_A^{(k)}\}^\monD)$$
$$
=1\cdot \left({1\over 2}+\sum_{i\in A} k\cdot {1\over 2nk} \right)-\sum_{i\in A^c} {1\over 2nk}
={1\over 2}+{k|A|\over 2nk}-{n-|A|\over 2nk}={n(k-1)+|A|(k+1)\over 2nk},$$
since $ \PP(\e_{A}^{(k)}+s_i,\{\e_A^{(k)}\}^\monD)= \PP(\e_{A}^{(k)}+s_i,\e_A^{(k)})$.
\smallskip\par
\noindent
Now, for the probability of transition from $\e_{A}^{(k)}$ to $\e_{A\cup\{j\}}^{(k)}, j\notin A$ we obtain
$$\PP^*(\e_{A}^{(k)},\e_{A\cup\{j\}}^{(k)})={H(\e_{A\cup\{j\}}^{(k)})\over H(\e_A^{(k)})}  \sum_{\e\succeq \e_{A\cup\{j\}}^{(k)}} 
\mu(\e_{A\cup\{j\}}^{(k)},\e)\PP(\e, \{\e_A^{(k)}\}^\monD).$$
The only state $\e$ for which $\PP(\e, \{\e_A^{(k)}\}^\monD)>0$  is $\e=\e_{A\cup\{j\}}^{(k)}$, thus (using (\ref{eq:HeAk})) we have

$$\PP^*(\e_{A}^{(k)},\e_{A\cup\{j\}}^{(k)})=(k+1)  \mu(\e_{A\cup\{j\}}^{(k)},\e_{A\cup\{j\}}^{(k)})\PP(\e_{A\cup\{j\}}^{(k)}, \{\e_A^{(k)}\}^\monD)={k+1\over 2nk}.$$
%$$={H(\e_{B_j}^{(k)})\over H(\e_A^{(k)})}\cdot 1\cdot {1\over 2nk}
%={(k+1)^{|B_j|}\over (k+1)^{|A|}}\cdot {1\over 2nk}
%={(k+1)^{|A|+1}\over (k+1)^{|A|}}\cdot {1\over 2nk}={k+1\over 2nk}$$
This completes our argument since all other transitions  have probability 0, which is clear from the following summation
$$\PP^*(\e_A^{(k)},\e_A^{(k)})+\sum_{j\in A^c}\PP^*(\e_A^{(k)},\e_{A\cup\{j\}}^{(k)})=
{n(k-1)+|A|(k+1)\over 2nk}+(n-|A|)\cdot {k+1\over 2nk}$$
$$={ n(k-1)+n(k+1)+|A|(k+1)-|A|(k+1)\over 2nk}=1.$$
Note that  the dual chain  starts at the minimal state which is also of the form $\e_A^{(k)},$ namely with $A=\emptyset$.